\documentclass{article}

\usepackage{PRIMEarxiv}

\usepackage[utf8]{inputenc} 
\usepackage[T1]{fontenc}    
\usepackage{hyperref}       
\usepackage{url}            
\usepackage{booktabs}       
\usepackage{amsfonts}       
\usepackage{nicefrac}       
\usepackage{microtype}      
\usepackage{lipsum}
\usepackage{fancyhdr}       
\usepackage{graphicx}       
\graphicspath{{media/}}     

\newtheorem{theorem}{Theorem}[section]

\usepackage{listings}
\usepackage{amsmath}
\usepackage{verbatim}
\usepackage{amssymb}
\usepackage{pdfrender}
\DeclareRobustCommand*{\pmbb}[1]{%
  \textpdfrender{
    TextRenderingMode=Stroke,
    LineWidth=.1pt,
  }{#1}%
}
\DeclareUnicodeCharacter{2212}{-}
\usepackage{float}

\title{ On the Boundedness solutions  of the difference equation $x_{n+1}=a x^\alpha_{n}+bx^\alpha_{n-1},0<\alpha \leq2$ and its application in  medicine
}

\author{
  Zeraoulia Rafik \\
  University of batna2.Algeria\\
  Departement of mathematics \\
  yabous ,khenchela\\
  \texttt{\{Author1\}r.zeraoulia@univ-batna2.dz} \\
   \And
  Alvaro humberto Salas \\
  Universidad national de Colombia \\
  departement of physics\\
  Bogota Colombia\\
  \texttt{ahsalass@unal.edu.co} \\
   \And
  Lorenzo Martinez\\
  Universidad national de Colombia \\
  departement of mathematics\\
  Manizales,Caldas\\
  \texttt{ljmartinezhe@unal.edu.co} \\
}

\begin{document}
\maketitle

\begin{abstract}
Recently ,mathematicians have been interested in studying the theory of discrete dynamical system, specifically   difference equation, such that  considerable works about  discussing  the behavior properties  of its solutions (boundedness and unboundedness)  are discussed and published in many areas of mathematics which involves several interesting results and applications in applied mathematics and physics ,One of the most important discrete dynamics which is become of interest for researchers in the field is the rational dynamical system .In this paper we may discuss qualitative behavior and properties of the difference equation 
$x_{n+1}=ax^2_{n}+bx^2_{n-1}$ with $a$ and $b$ are  two parameters and we shall show its application to medicine
\end{abstract}


\section{Introduction}
The theory of difference equations finds many applications
in almost all areas of natural science \cite{Dan:06}. 
increasingly clearly
emerges the fundamental role that difference equations with
discrete and continuous argument is played for understanding
nonlinear dynamics and phenomena also it is used for combinatorics and in the approximation of solutions of partial differential equations \cite{far:17}. The increased interest in difference equations is partly
due to their ease of handling. A minimum is enough
computing and graphical tools to see how the
solution of difference equations trace their bifurcations with
changing parameters \cite{Josef:08}. Thus opens a complex understanding as well invariant manifolds for linear and nonlinear dynamical systems.nonlinear difference equations and systems are of wide interest due to their applications in real life. Such equations appear naturally as the mathematical models which describe biological, physical and economical phenomena.

Although difference equations have very simple forms, however, it is extremely difficult to understand completely the global behavior of their solutions.y the global behaviors of their solutions. One can refer to \cite{RP:92},\cite{G.ladas:08},\cite{E.A:05} and
the references therein. Difference equations have always played an important role in the
construction and analysis of mathematical models of biology, ecology, physics, economic
processes, etc. The study of nonlinear rational difference equations of higher order is of
paramount importance, since we still know so little about such equations.
In this paper \cite{A khan:20} A. Q. Khan and S. M. Qureshi discussed Dynamical properties of some rational systems of difference equations such as they  explored the equilibrium points, local and global dynamics, rate of convergence, instability and boundedness of positive solution of some rational systems of difference equations .As application to modern science ,namely, in mathematical biology they  also explored the local dynamics about equilibrium points of the discrete-time Levin's model .In the meanwhile A. Q. Khan studied and discussed Global dynamics of a $3 \times 6$ exponential system of difference equations defined as:

$$\begin{cases} x_{n+1}=\frac{\alpha_{1}+\beta_{1}\exp(-x_n)}{\gamma_1+y_{n-1}} \\  y_{n+1}=\frac{\alpha_{2}+\beta_{2}\exp(-y_n)}{\gamma_2+z_{n-1}}  \\ z_{n+1}=\frac{\alpha_{3}+\beta_{3}\exp(-z_n)}{\gamma_3+y_{n-1}} \end{cases}$$
$n=0,1,....$

where parameters  $\alpha_i,\beta_i,\gamma_i (i=1,2,3)$ and initial conditions $x_i,y_i,z_i (i=0,−1)$ are nonnegative real numbers. 

In \cite{R.Abo:2017} .R. Abo Zeid has discussed the  global behavior of all solutions of the difference equation:
\begin{equation}\label{eq:1}
    x_{n+1}=\frac{x_n x_{n-1}}{ax_n+bx_{n-1}},n\in \mathbb{N}_{\pmbb{0}}
\end{equation}
where $a, b$ are real numbers and the initial conditions $x_{−1}, x_{0}$ are real numbers .In this paper, we discuss the global behavior of the difference equation :
\begin{equation}\label{eq:2}
    x_{n+1}=A x_n ^\alpha +B x_{n-1} ^{\alpha}
\end{equation}

with :$A , B$ are two parameters real numbers and $\alpha$ is a real number such that :

$0<\alpha \leq 2$, For $\alpha =1$ the dynamics defined in (\ref{eq:2}) is deduced from (\ref{eq:1}) by the substitution $y=\frac{1}{x_n}$ , the  global behavior of all solutions of (\ref{eq:1}) is discussed as well in  \cite{R.Abo:2017} by R. Abo-Zeid ,The difference equation (\ref{eq:2}) for $\alpha=1$ ,namely , 
\begin{equation}\label{eq:3}
    y_{n+1}=b y_n+ay_{n-1} ,n\in \mathbb{N}
\end{equation}
The characteristic equation of equation (\ref{eq:3}) is :\begin{equation}\label{equ:4}
\lambda^2 -b \lambda -a=0
\end{equation}
Equation (\ref{equ:4}) has two roots : $\lambda_1=\frac{b-\sqrt{b^2-4ac}}{2a},\lambda_2=\frac{b+\sqrt{b^2-4ac}}{2a}$
The form of the solution should be according to the value of the quantity $b^2+a^2$,The following theorem \cite{Elaydi:05} is useful in studying the solutions of the difference equation (\ref{eq:3}) 
\begin{theorem}
The following statements holds:
\begin{itemize}
    \item 1) All solutions of (\ref{eq:3}) oscillates (about zero) if and only if the characteristic equation has no positive roots.

\item 2) 
All solutions of (\ref{eq:3}) converge to zero if and only if $\max \{\lambda_1,\lambda_2 ,\}<1
$

\end{itemize}
\end{theorem}
For boundedness solutions of (\ref{eq:3})  one can refer to \cite{R.Abo:2017} .Now for $\alpha=2$ which it is the aim of this paper we are ready to do some analysis and discussion about the global behavior of solutions of this difference equation:
\begin{equation}\label{eq:5}
        y_{n+1}=b y_n ^2+ay_{n-1}^2 ,n\in \mathbb{N}
\end{equation}

\section{Analysis and discussion:}

\textbf{Case:1}
$|a|,|b|< 1$ ,for this case we may use a nice trick just we assume $a, b$ are two bounded functions such that : $a=\sin(\theta) ,b=\cos (\theta),\theta \in \mathbb{R}$ then the dynamics defined in (\ref{eq:5}) become :
\begin{equation}\label{eq:6}
y_{n+1}=\cos (t) y_n ^2+\sin (t) y_{n-1}^2 ,n\in \mathbb{N},t \in \mathbb{R}
\end{equation}
Now ,we may ask for which values of $\theta$ does this equation $x_{n+1}=\cos(\theta)x^2_{n}+\sin(\theta)x^2_{n-1}$ have bounded solutions?.

\begin{figure}[H]
    \centering
    \includegraphics[width=0.4\textwidth]{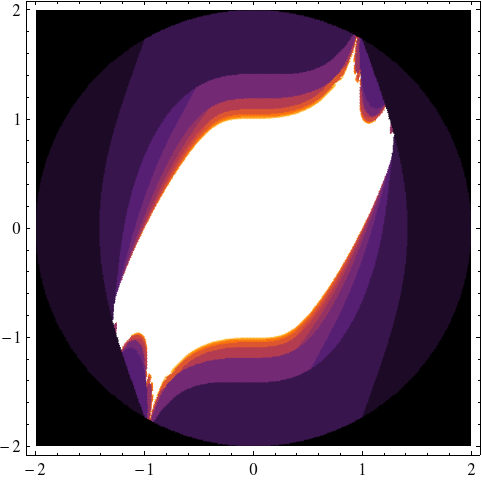}
    \caption{Bounded solution for $x_{-1}=0,r=x_{0}$}
    \label{fig:Nice}
\end{figure}

We  ran a small computation: The below plot,see figure (\ref{fig:Nice}) is created as follows:
For each point $r(\cos \theta, \sin \theta)$, We use $x_{-1} = 0$, $x_{0}= r$
as initial values, and $\theta$ as the parameter.

The white area is where the iterations is less than $2$ (bounded), the first 30 iterations. As we see, a small circle around the origin is white, meaning that there are small initial values that is bounded for every $\theta$. Now, this is not a proof, but the picture suggests this.

Changing the cut-off to $40$ iterations does not change the picture much.

For analytical proof we have :Let $r=1/\sqrt{2}$. If $x_{-1},x_{0} \leq r$, then the sequence is bounded for all $\theta$.

We note that: $|\cos(\theta)x_{n-1}^2 + \sin(\theta)x_{n}^2| \leq r^2(|\cos \theta|+|\sin \theta|) \leq r^2 \sqrt{2} \leq 1/\sqrt{2}$, and the statement follows from induction.

\textbf{Case:2}

\begin{figure}[H]
    \centering
    \includegraphics[width=0.4\textwidth]{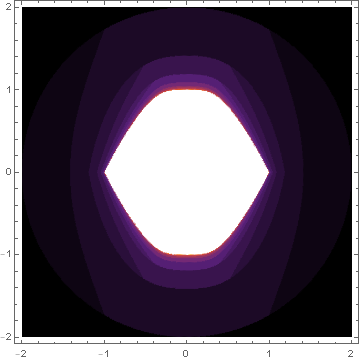}
    \caption{new intial values that bounded solutions for every $\theta ,x_{-1}=0,r=x_0$}
    \label{fig:G}
\end{figure}

$|a|,|b|= 1$, for the case $|a|=|b|=1$ we note after runing with the same initial conditions for each point $r(\cos \theta, \sin \theta)$  (iteration of the plot (20 iteration)  )  that the obtained  picture changed ( Rotation of picture) and the white region become larger than it for the precedent case ,namely, \textbf{Case:1} this mean that we have another new initial values  that is bounded for every $\theta$.(The number of initial values is increasing).see figure \ref{fig:G}

\textbf{Case3}
In this case we may assume :$a'=|a|k ,b'=|b|k,k\in \mathbb{R}$ and $|a'|, |b'| >1$, then our dynamics become :
\begin{equation}\label{eq:7}
y_{n+1}=|\cos (t)| y_n ^2+|\sin (t)| y_{n-1}^2 ,n\in \mathbb{N},t \in \mathbb{R}
\end{equation}
With the same data and conditions that we have used in both of case1 and case 2, we noticed after runing using mathematica  the white region become smaller  whenever $k$ being larger this indicate that the dyanmics defined in (\ref{eq:7}) would be unbounded,for numerical evidence we just take , $k=5,20,-15$  as shown in below plots (figure (3+4+5)).
\begin{figure}[H]
    \centering
    \includegraphics[width=0.4\textwidth]{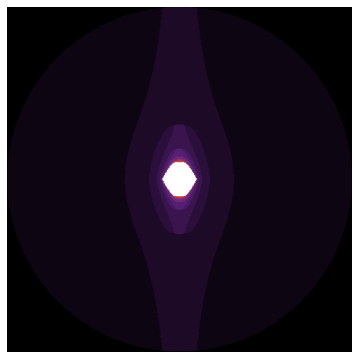}
    \caption{boundedness solutions for the dynamic (\ref{eq:7}),$k=5$ }
    \label{fig:A}
\end{figure}

\begin{figure}[H]
    \centering
 \includegraphics[width=0.4\textwidth]{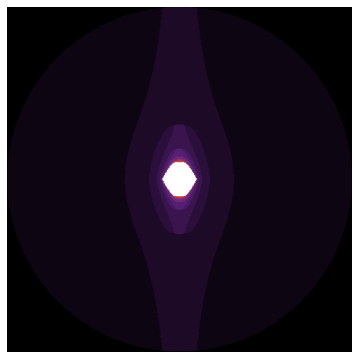}
    \caption{boundedness solutions for the dynamic (\ref{eq:7}),$k=15$ }
    \label{fig:J}
\end{figure}

\begin{figure}[H]
    \centering
 \includegraphics[width=0.4\textwidth]{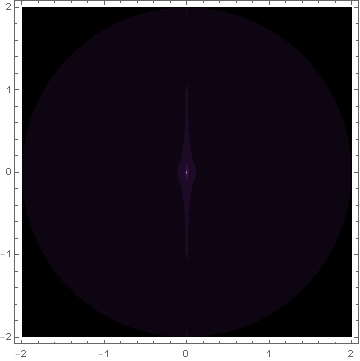}
    \caption{boundedness solutions for the dynamic (\ref{eq:7}),$k=-150$ }
    \label{fig:K}
\end{figure}

We noticed For the remainder case when  $ 0<\alpha \leq1$ the iteration of the dynamics (\ref{eq:2}) show us  much  intinial values that is bounded for every $\theta$ , we took $\alpha =\frac12,\frac13,\alpha=1$ as comparative examples as shown in below figures: 

\begin{figure}[H]
    \centering
    \includegraphics[width=0.2\textwidth]{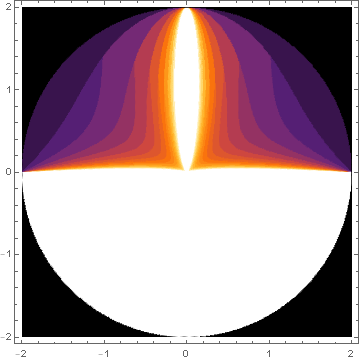}
    \caption{bounded solutions for $\alpha=1$ (white circle arround origin)}
    \label{fig:R}
\end{figure}

\begin{figure}[H]
    \centering
    \includegraphics[width=0.4\textwidth]{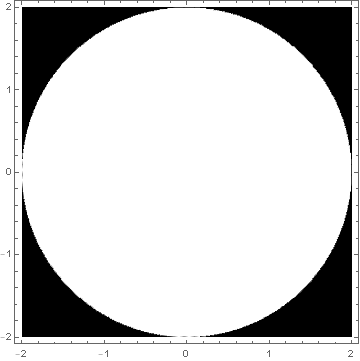}
  \includegraphics[width=0.4\textwidth]{293075314_531845868694694_5419669993824061089_n.png}
    \caption{bounded solution for $\alpha=\frac{1}{2},\frac{1}{3},\cdots$ (The circle become white at all this indicate all solution of our dynamics are bounded for any arbitary initial values )}
    \label{fig:Rafik}
\end{figure}

\section{Stability and analysis}

Dynamical modeling is the study of change and changes take place everywhere in life. As a result dynamical systems have a wide range of application areas in applied science and engineering. With these systems, real life situations can be turned into the language of mathematics. If one can create a good model for a real life situation, we will be able to predict the future states of the system by simply iterations according to this model.

Stability, like one of the most important concepts in Discrete Dynamical Theory, can tell much about the behavior of the dynamical system. In this section a symbolic Mathematica package for analysis and control of chaos in discrete  two  dimensional dynamical nonlinear systems, is presented. There are constructed some computer codes to find stability types of the fixed points, covering the stability of the one-dimensional nonlinear dynamical systems. Applications are taken from chemical kinetics and population dynamics (logistic model).Since our dynamics is two dimonsional discret dynamical system it is enough to run the below  code to get satibility and oscillation points arround origing...

\textbf{example one}
Let us apply that code for our dynamics:
let $a=2, b=9$

\begin{verbatim}
    twoDimStab[2 x^2 + y + 1, 9x^2 - 1]

       1     112
{Null, {--, -(---)} -> oscillatory source}
        11    121
\end{verbatim}

One can do more examples just applying the above code with change of variable $a$ and $b$.

For the periodicity of our dynamic for $\alpha=2$ it is already discussed in  \cite{R.Abo:2017},The dynamics is two period solution

\section{PHASE PLANE DIAGRAMS}

Continuous systems are often approximated as discrete processes, meaning that we look only at the solutions for positive integer inputs. Difference equations are recurrence relations, and first order difference
equations only depend on the previous value. Using difference equations, we can model discrete dynamical systems. The observations we
can determine, by analyzing phase plane diagrams of difference equations, are if we are modeling decay or growth, convergence, stability,
and equilibrium points.In this section we may give a diagram plot of our dynamics for some values $a$ and $b$ and $\alpha=2$ , The dynamics defined in (\ref{eq:2}) can be written as sytem of two difference equations :
\begin{equation}\label{S}
\begin{cases}
 x(t+1)=ax_t^2+y(t)+1\\
y(t+1)=b x_t^2-1
\end{cases}
\end{equation}

The phase portrait plot of the dynamic (\ref{S}) is shown in  below figures,
Here is the mathematica code one can change values of a,b to analyse the slop field 
\begin{figure}[H]
    \centering
    \includegraphics[width=0.5\textwidth]{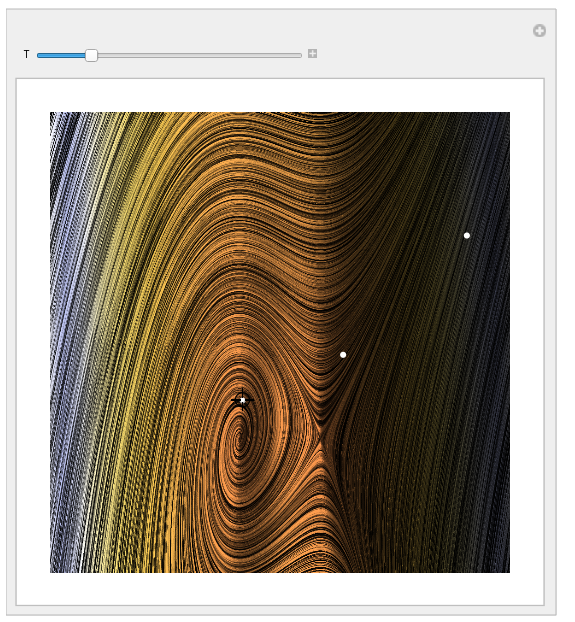}
    \caption{diagram plot of nonlinare dynamics (\ref{S}), for $a=0.5, b=2$}
    \label{fig:M}
\end{figure}

\section{Time series analysis and modeling (application to medicin)}

Nowadays mathematics is being successfully applied to a number of important fields in medicine including biofluids, cardiovascular diseases, clinical schedules and tests, data analysis, drug design and discovery, epidemiology, genetics, image processing, immunology, instrumentation, microbiology \cite{Small:17}, neuroscience, oncology, virology and more. The list of tools includes virtually the whole of applied mathematics. To cite the most familiar ones: difference equations and discrete-time dynamical systems, information and coding theory, graph and network theory, integral transforms, numerical and computational mathematics, ordinary differential equations and continuous-time dynamical systems, partial differential equations, stochastic and time-delay differential equations, statistics, probability and time-series analysis. All this research has contributed to and continues to increasingly contribute both to better understand medical phenomena and to finding practical ways of action.

Time series modeling for predictive purpose has been an active research area of machine learning
for many years \cite{Fat:19}. However, no sufficiently comprehensive and meanwhile substantive survey was
offered so far. This survey strives to meet this need. A unified presentation has been adopted for
entire parts of this compilation.
Time series data are amongst the most ubiquitous data types that capture information and record
activity in most aeras. In any domain involving temporal measurements via sensors, censuses, transaction
records, the capture of a sequence of observations indexed by time stamps first allows to provide insights
on the past evolution of some measurable quantity. Beyond this goal, the pervasiveness of time series has
generated an increasing demand for performing various tasks on time series data (visualization, discovery
of recurrent patterns, correlation discovery, classification, clustering, outlier detection, segmentation,
forecasting, data simulation).

Time series classification (TSC) is one of data minings persistent challenges.
Applications of TSC abound in fields including agriculture, medicine, and engine prognostics \cite{Abdoli:18},and \cite{Yu:21} a
common goal being to detect instances of sub optimal behavior or decreased health (biologically or
mechanically) as just one real-world example. Dozens of new TSC algorithms were introduced in the last
four years alone \cite{Bagnall:17}. This trend has been intensified by the increasing availability of real-world
datasets. In fact, the classification of any inherently ordered data (temporally or otherwise) can be
treated as a TSC problem \cite{Gam:17} , making for a vast breadth of real-world applications.
Deep learning methods have shown suitability for time series classification in the health and medical domain, with promising results for electrocardiogram data classification. Successful identification of myocardial infarction holds life saving potential and any meaningful improvement upon deep learning models in this area is of great interest.In this section we show that the discussed dynamics ,namely , (\ref{S})  present new classes of heartbeat ,namely, A new dynamical model for Generating Synthetic Electrocardiogram Signals,The electrocardiogram (ECG) is a time-varying signal reflecting the ionic current flow which causes the cardiac fibers to contract and subsequently relax. The surface ECG is obtained by recording the potential difference between two electrodes placed on the surface of the skin. A single normal cycle of the ECG represents the successive atrial depolarization/repolarization and ventricular depolarization/repolarization which occurs with every heartbeat. These can be approximately associated with the peaks and troughs of the ECG waveform labeled $P, Q, R, S$, and $T$ as shown in figure 9

\begin{figure}[H]
    \centering
    \includegraphics[width=0.5\textwidth]{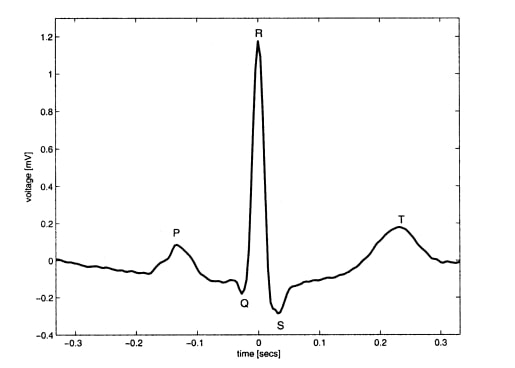}
    \caption{Morphology of a mean PQRST-complex of an ECG recorded from a
normal human}
    \label{fig:N}
\end{figure}
Extracting useful clinical information from the real (noisy) ECG requires reliable signal processing techniques  \cite{Goldberger:77}. These include R-peak detection \cite{Jpn:85}, \cite{Tai:95}, QT-interval detection 
and the derivation of heart rate and respiration rate from the
ECG \cite{P.E:03}. The RR-interval is the time between successive
R-peaks, the inverse of this time interval gives the instantaneous heart rate. A series of RR-intervals is known as a RR tachogram and variability of these RR-intervals reveals important information about the physiological state of the ECG signal. The ECG may be divided into the following sections.
\begin{itemize}
\item
 Q-wave: A small low-voltage deflection away from the
baseline caused by the depolarization of the atria prior to
atrial contraction as the activation (depolarization) wavefront propagates from the SA node through the atria.
\item
 PQ-interval: The time between the beginning of atrial
depolarization and the beginning of ventricular depolarization.

\item
 QRS-complex: The largest-amplitude portion of the ECG,
caused by currents generated when the ventricles depolarize prior to their contraction. Although atrial repolarization occurs before ventricular depolarization, the latter
waveform (i.e. the QRS-complex) is of much greater amplitude and atrial repolarization is therefore not seen on
the ECG.
\item
 QT-interval: The time between the onset of ventricular
depolarization and the end of ventricular repolarization.
Clinical studies have demonstrated that the QT-interval
increases linearly as the RR-interval increases . Prolonged QT-interval may be associated with delayed ventricular repolarization which may cause ventricular tachyarrhythmias leading to sudden cardiac death 
\item
ST-interval: The time between the end of S-wave and the
beginning of T-wave. Significantly elevated or depressed
amplitudes away from the baseline are often associated
with cardiac illness.
\item
 T-wave: Ventricular repolarization, whereby the cardiac
muscle is prepared for the next cycle of the ECG.
\end{itemize}
Analysis of variations in the instantaneous heart rate time
series using the beat-to-beat RR-intervals (the RR tachogram)
is known as HRV analysis \cite{AJ:95}, \cite{North:96}. HRV analysis has been shown to provide an assessment of cardiovascular disease \cite{M H Crawford:96}.A dynamical model based on three coupled ordinary differential equations is introduced which is capable of generating realistic synthetic electrocardiogram (ECG) signals.The dynamical equations of motion are given by a set of three ordinary differential equations defined as the following :
\begin{equation}\label{N1}
\begin{cases}
 \dot{x}=\alpha x-\omega y\\
 \dot{y}=\alpha y+\omega x\\
 \dot{z}=-\sum_{i\in \{P,Q,R,S,T \}}\exp \bigg(-\frac{-\Delta \theta_i^2}{2b_i^2}\bigg)-(z-z_0)
\end{cases}
\end{equation}
where $\alpha=1-\sqrt{x^2+y^2},\Delta \theta_i=(\theta-\theta_i)\bmod 2\pi,\theta=\arctan 2(y,x)$,  (the four quadrant arctangent of the real parts of elements of $x$ and $y$ , with  $-\pi \leq \arctan 2(y,x)\leq \pi$) and $\omega$  is the angular velocity of the trajectory as it moves around the limit cycle. Baseline wander was introduced by coupling the baseline value $z_0$ in (\ref{N1}) to the respiratory frequency $f_2$ using :$z_0(t)=A\sin (2\pi f_2 t)$ where $A=0.15mV$

\begin{figure}[H]
    \centering
    \includegraphics[width=0.5\textwidth]{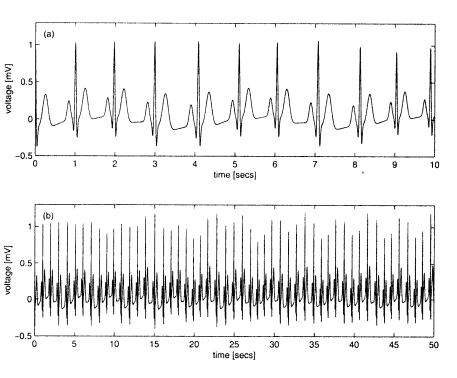}
    \caption{ ECG generated by dynamical model: (a) 10 s and (b) 50 s}
    \label{fig:Np}
\end{figure}
\begin{figure}[H]
    \centering
    \includegraphics[width=0.5\textwidth]{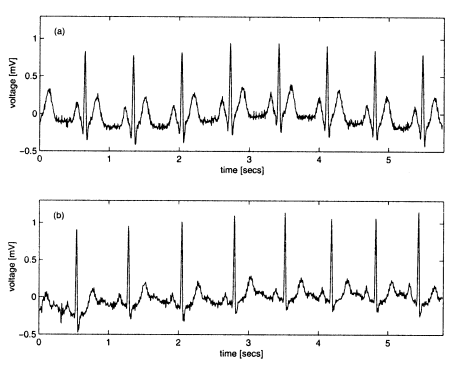}
    \caption{Comparison between (a) synthetic ECG with additive normally
distributed measurement errors and (b) real ECG signal from a normal human.}
    \label{fig:Npd}
\end{figure}

\section{Introduction to Electrocardiography}

An electrocardiogram is a recording of the electrical activity of the heart.The heart can be viewed as a three-dimensional vector. Therefore, its electrical activity can, in theory, be recorded by three orthogonal leads. In practice, however, a standard clinical EKG is recorded with 12 leads: 6 limb leads and 6 precordical leads.A normal EKG reflects the electrical activity of each of the four heart chambers: left and right ventricles and left and right atria (Fig. 1).  

\begin{figure}[H]
    \centering
    \includegraphics[width=1\textwidth]{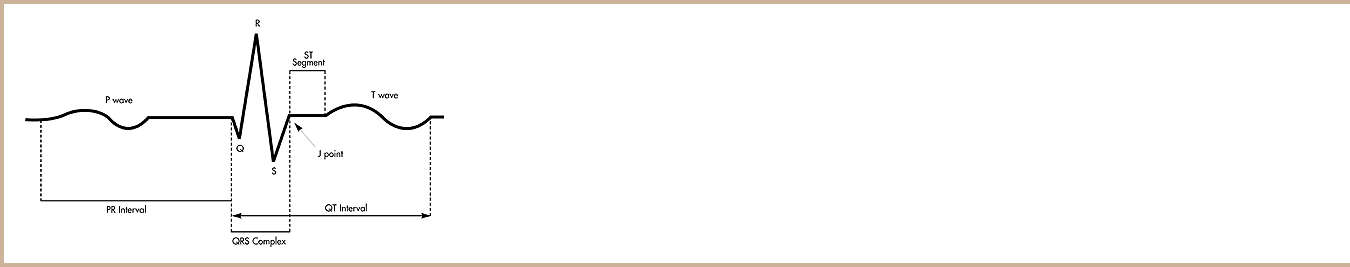}
    \caption{A normal EKG reflects the electrical activity of each of the four heart chambers: left and right ventricles and left and right atria }
    \label{fig:Algeria}
\end{figure}
The P wave marks the beginning of atrial depolarization. The onset of the Q wave is produced by the stimulus from the Purkinje system. The depolarization of the ventricle produces the R wave. The S wave is initiated when the excitation reaches the base of the ventricle. The T wave is produced by ventricular repolarization (Tompkins 1993; Goldman 1992).
The synthetic ECG (Figure 10) illustrates the modulation of the
QRS-complex due to RSA and Mayer waves. Observational
uncertainty is incorporated by adding normally distributed
measurement errors with mean zero and standard deviation
$0.025 mV$ [Figure 11], yielding a similar signal to a segment of
real ECG from a normal human [Figure 11].
Electrocardiogram (ECG) detection is currently the most
effective and direct way to detect ECG signals \cite{Z.Tang}. At present,
the diagnosis of cardiac diseases is mainly determined by medical doctors and clinicians through manual detection and ECG
analysis.ECG is a diagnostic technology that records the electrocardiography activities of the heart in a certain time unit
through the chest of biological objects. It collects and records
the electrodes connected to the skin of specific parts of biological objects and preserves the relevant contents in a certain form  \cite{Wolf:97}.

The pre-ejection-period (PEP) is the time span between the depolarization of the left ventricle (R onset) and opening of the aortic valve (B point). The R onset is the beginning of the Q-wave signal; it indicates the beginning of the depolarization and can be picked up from the ECG signal. As signature for the beginning of the Q wave, we take the minimum of the ECG's second derivative. Its peak indicates the maximum curvature at the transition of the ECG signal into the Q wave. However, as the Q wave is relatively small, other signatures in the ECG signal can be misinterpreted as the R onset in an automated evaluation of the data. In general, first and higher-order derivatives of noisy signals suffer from containing spurious peaks. We therefore restrict the possible occurrences of R onset to a time window after the, easily identifiable, R peak (peak of the QRS complex). Within that window, the R onset is typically seen as a clear negative peak of the ECG signal's second derivative, which can be located reliably and with high precision and thus allows a reliable identification of the Q wave onset. Heart rate (HR) is then calculated from the time difference of subsequent R points.\cite{Ari:15}

The time point for the opening of the aortic valve (B point) is derived from the impedance cardiogram (ICG). The impedance Z, and thus the ICG, is sensitive to a variation of blood volume in the thorax. The first derivative, $\frac{dZ}{dt}$, corresponds to blood flow. The second derivative  $\frac{d^2 Z}{dt^2}$, in turn, corresponds to a change of the blood flow and is thus indicative for the opening of the heart valves. The B point is the onset of the aortic valve's opening, indicated by a negative peak in the third derivative, $\frac{d^3 Z}{dt^3}$. While, compared to ECG, the ICG signal is smooth and devoid of characteristic spikes, its first, second, and third derivative show distinct features. As selection criterion for picking the correct peak of the third derivative, we use the, easily identifiable, peak of the first derivative, $\frac{dZ}{dt}$. The B point is obtained as the minimum of $\frac{d^3 Z}{dt^3}$ that occurs just before the maximum in $\frac{dZ}{dt}$. This strategy allows for an automated evaluation of the PEP interval for the large data sets, with few outliers and the required precision.
To calculate the derivatives of the measured signals, we use the Savitzky-Golay filter \cite{Jianwen LuoKui:05}. This method allows data smoothing, while keeping intact signatures like peaks, and the simultaneous determination of derivatives. Similar to a moving average, a moving section of the data is selected. However, instead of a simple averaging, the algorithm fits a polynomial of given degree to the selected sequence. Then, one point of the fitted polynomial (usually the central point) is taken as value for the smoothed curve. Higher derivatives are taken from the corresponding derivatives of the fitted polynomial at the respective point. The Savitzky-Golay filter is implemented numerically by a list convolution with a kernel. That kernel is calculated in advance for the number of points for the moving fitting, the order of the polynomial, and the order of the derivative. We use a kernel length of 100 points, corresponding to a time interval of 50 ms, and a 3rd-order polynomial for all kernels and for the ICG and ECG signals. The third derivative of the ICG signal is calculated from the first derivative of the ICG signal, which, together with Z, is provided by the Biopac MP36 system (i.e., the system we used to measure ICG/ECG). We ensure that no time lag gets introduced between the ICG and ECG signals and their derivatives by the Savitzky-Golay filter.  Thus, PEP and LVET data get extracted from the ICG and ECG measurements in a semi-automated way and with a by-heartbeat resolution. The Mathematica Notebook output is stored in a text file, with every row containing a timestamp, the corresponding length of cardiac PEP, LVET, and HR for each heartbeat. Here is the Graphical display of $Z,\frac{dZ}{dt}$,EKG 

\begin{figure}[H]
    \centering
    \includegraphics[width=0.5\textwidth]{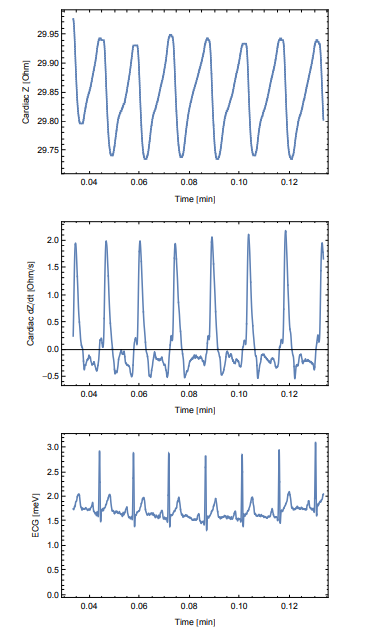}
    \caption{Graphical display of $Z,\frac{dZ}{dt}$ ,EKG , Istart=1000,Iend =16000}
    \label{fig:Npdt}
\end{figure}
\section{Medical interpretation of our discret dynamics }
To analyze ECG signal, the most necessary step is to extract its QRS wave group. The QRS complex reflects changes in depolarization potential and time of the left and right ventricles. Considering the robustness and stability, For analysis of ECG signal using delay differential equation (DDE's) is computationally fast one can refer to the system of ODE in (\ref{N1}), and could be the basis for a real time diagnostic system. DDEs reveal non-linear properties as well as spectral properties of the data as shown in this paper \cite{Claudia Lainscsek:13} by Claudia Lainscsek and Terrence J. Sejnowski they analyzed ECG signal using delay differential equation which leads to the good classification of Electrocardiogram ,may that classification wouldn't work as well using discret map ,some authors discussed Electrocardiogram Signal Classification in the Diagnosis of Heart Disease Based on RBF Neural Network \cite{Yan Fang:22} such that they extracted QRS wave using discret map (they used difference equations to analyze ECG ),Recall that in our paper we are interested to the behavior of the following dynamics :$$x_{n+1}=ax^2_{n}-bx^2_{n-1} ,n=0,1,....$$ ,$a,b$ two real parameters ,It is known that ECG signal has a very obvious periodicity, Taking $T$ as sampling period we may rewrite our dynamics using $T$ as :
\begin{equation}\label{Raf}
    x((n+1)T)=ax^2 (nT)-bx^2(T(n-1),n=0,1,...
\end{equation}
The high-frequency characteristics can be enhanced by a nonlinear square function \cite{Yan Fang:22},see page 3 equation 4, whose equation can be expressed as :
\begin{equation}\label{Rafik}
    y((nT)=x^2 (nT)
\end{equation}
This means strongly that our discret  dynamics which is defined in (\ref{Raf}) can be interpreted in the medical  point of view as :Enhanced high frequencies  regarding  the behavior of the heartbeat in  short time, we may consider  the discret time defined in (\ref{Raf}) as a discretized boundary value problem (delay differential equation) using some standrad numerical methods like finite difference  method  in 1D (dimension)  , we may give a short analysis of the ECG signal which is produced by the dynamics defined in (\ref{Raf}), Assume the correspending non linear boundary value problem which we wanted to solve satisfies  $x_1 =0$ and $x_{n+1} = 1$  using the Finit difference technique \cite{R.P.Ag:85}. The first step is to partition the domain $[0,1]$ into a number of sub-domains or intervals of length $h$. So, if the number of intervals is equal to $n$, then $nh = 1$. We denote by $x_i$ the interval end points or nodes, . In general, we have $x_i = (i-1)h$, $i=1,2,\cdots n+1$. Let us denote the concentration at the ith node by $C_i$, for short enough time (say between $t_i, t_{i+1}$) the dynamics (\ref{Raf}) become as a simple linear differential equation such that :
for $h\to 0$ we have $x((n+1)T) \to \dot{x} \leq 0$ implies that we have a decreasing frequencies coming up to a constant signal $x(t)$ thus ,the dynamic (\ref{Raf}) can be interpreted as heart attack or heart failure  ,In general that case indicate  Cardiac diseases ,(see figure in example 1), Enhanced high frequency appear whenever $\dot{x} > 0$ which means derivative of the analyzed signal is always positive (one can refer to figures in examples (2+3+4)).It is hard to analyze ECG signal using coupled signal which is defined in (\ref{S1})

we have plotted many figures with  some values of $a$ and $b$ such that we call for medical interpretation  $a$ is a factor of life  (Electrocardiogram (EKG)) ,identification of that factor(EKG)  present attempt to improve patient survival and $b$ may present blood loss which causes the cardiac arrest.we may introduce a new parameter here $\sigma$ as a toxine factor  this  for good modelization of the discussed phenomena.
Now,let us try showing heartbeat classes  ,The dynamics (\ref{S}) for $\alpha=2$ can be written as :
\begin{equation}\label{S1}
\begin{cases}
 x(t+1)=ax_t^2+\sigma y(t)+1\\
y(t+1)=b x_t^2-1
\end{cases}
\end{equation}

We have noticed that For heartbeat classes plots  we should have  : $a \leq 0.15$,we may try playing with $b$ and $\sigma <0 $ (should be negative) and for $x_0=0.07,y_0=0.08$ initial values  ,increasing values of $\sigma $ means attempt to reduce toxine factor and increasing values of $b$ means raise factor of EKG thus a good attempt to improve patient survival ,for the reverse assumption we would have cardiac arrest record.
\\

\textbf{Example1}

we start with the first heartplot  , $a=0.15,\sigma =-0.45,b=0.45$

\begin{figure}[H]
    \centering
    \includegraphics[width=0.5\textwidth]{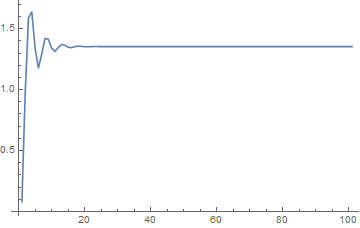}
    \caption{heartbeat classe for  $a=0.15,\sigma =-0.45,b=0.45,x_0=0.07,y_0=0.08 $}
    \label{fig:1}
\end{figure}

\textbf{Example2}
let :$a=0.15,\sigma =-0.45,b=0.45$
\begin{figure}[H]
    \centering
    \includegraphics[width=0.5\textwidth]{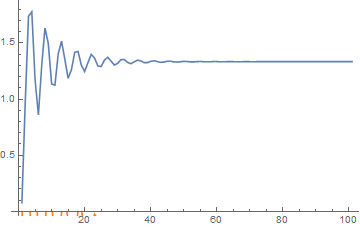}
    \caption{heartbeat case for  $a=0.15,\sigma =-0.6,b=0.5,x_0=0.07,y_0=0.08$}
    \label{fig:2}
\end{figure}
\textbf{Example3}
let :$a=0.15,\sigma =-0.65,b=0.58$
\begin{figure}[H]
    \centering
    \includegraphics[width=0.5\textwidth]{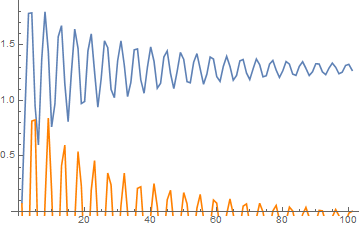}
    \caption{heartbeat case for  $a=0.15,\sigma =-0.65,b=0.58,x_0=0.07,y_0=0.08$}
    \label{fig:3}
\end{figure}

\textbf{Example4}
let :$a=0.15,\sigma =-0.65,b=0.57$
\begin{figure}[H]
    \centering
    \includegraphics[width=0.5\textwidth]{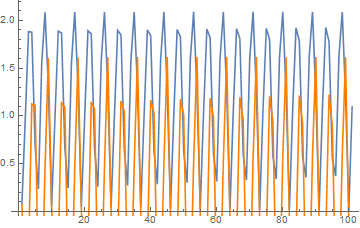}
    \caption{heartbeat case for  $a=0.15,\sigma =-0.75,b=0.6,x_0=0.07,y_0=0.08$}
    \label{fig:4}
\end{figure}

\section{Conclusion:}
This study of difference equations with public health applications to science develops the methodology for the solution of the general kth order linear difference equation using the generating function approach and computers tools like mathematica. It includes an examination of the dynamics of disease spread and containment in populations using illness-death models ,Cardiovascular disease is one of the major hazard to human health today. ECG stands for electrocardiogram and it is an important way for clinical diagnosis cardiovascular disease. The ECG
refers to the small voltages (~1mv) found on the skin as a result of electrical activity of the heart. 
These
electrical actions trigger various electrical and muscular activity in the heart. The health and function of
the heart can be measured by the shape of the ECG waveform. Typical heart problems are leaking valves
and blocked coronary arteries.in our paper we have  discussed a new discret dynamics which lead to discover a new illness-death model(Enhanced high frequencies model ) using time series diagram and ECG analyses with some data which are defined by iterative discret dynamics such as  heartbeat classes ,in particulary attempts to improve patient for being survival,we may call it model of life for more time.

\section{Data Availability
}

The data that support the findings of this study can be obtained from the corresponding author upon reasonable request.

\section{Conflicts of Interest
}
The author declare that they have no conflicts of interest.

\section{Authors’ Contributions}
All authors contributed to the study conception and design of that research . Material preparation, data collection,  analysis and medical interpretation  were performed by Zeraoulia Rafik  . The first draft of the manuscript was written by Zeraoulia Rafik  and all authors commented on previous versions of the manuscript. All authors read and approved the revised  manuscript.

\bibliographystyle{elsarticle-harv}

\begin{thebibliography}{7}
\expandafter\ifx\csname natexlab\endcsname\relax\def\natexlab#1{#1}\fi
\expandafter\ifx\csname url\endcsname\relax
  \def\url#1{\texttt{#1}}\fi
\expandafter\ifx\csname doi\endcsname\relax
  \def\doi#1{\texttt{#1}}\fi
\expandafter\ifx\csname urlprefix\endcsname\relax\def\urlprefix{URL: }\fi
\expandafter\ifx\csname doiprefix\endcsname\relax\def\doiprefix{DOI: }\fi

\bibitem[{Daniel J. Duffy(2006)}]{Dan:06}
Daniel J. Duffy, 2006.Finite Difference Methods in Financial Engineering: A Partial Differential Equation Approach, ISBN:=9780470858820, 9780470858820



\bibitem[{Y. Ordokhan1 ,S. Davaei far(2017)
}]{far:17}
Able, B., 1956.Approximate Solutions of Differential Equations by Using the Bernstein Polynomials
\newline\doiprefix\doi{10.5402/2011/787694}

\bibitem[{Josef Diblik
,Miroslava Ruzickova
, Barbora Vaclavikova
(2008)}]{Josef:08}
Josef Diblik
,Miroslava Ruzickova
, Barbora Vaclavikova, 2008. Bounded Solutions of Dynamic Equations on
Time Scales,{International Journal of Difference Equations (IJDE).
,ISSN 0973-6069 Volume 3 Number 1 (2008), pp. 61–69
}

\bibitem[{Charlie y Routh(1992)}]{RP:92}
Agarwal,Difference  equations  and  inequalities,1992.first edition, Marcel Dekker

 \bibitem[ Camouzis and G.Ladas(2008)]{G.ladas:08}
 Dynamics of Third-Order Rational Difference Equations;With Open Problems and Conjectures,  2008. Chapman and Hall/HRC Boca Raton,


\bibitem[{E.A.Grove and G. Ladas(2005)}]{E.A:05}
E.A.   Grove and G. Ladas,Periodicities in Nonlinear Difference Equations,2005. Chapmanand Hall/CRC, 2005

\bibitem[R. Abo-Zeid
(2017)]{R.Abo:2017}
R. Abo-Zeid
On the solutions of a second order
difference equation,2017.{Mathematica Moravica
Vol. 21, No. 2 (61-73}
 
 
 \bibitem[S.Elaydi(2005)]{Elaydi:05}
  S.Elaydi, An Introduction to Difference Equations,2005. {Third Edition, Springer, New
York}.


\bibitem[Fatoumata Dama and Christine Sinoquet
(2019)]{Fat:19}
Fatoumata Dama and Christine Sinoquet,2019.
Time Series Analysis and Modeling to Forecast: a Survey .{LS2N / UMR CNRS 6004, Nantes University, France}



\bibitem[Abdoli, A., Murillo, A. C., Yeh, C. C. M., Gerry, A. C., and  Keogh, E. J. (2018, December)]{Abdoli:18}
Abdoli, A., Murillo, A. C., Yeh, C. C. M., Gerry, A. C.,  Keogh, E. J. 2018. Time series
classification to improve poultry welfare. {In 2018 17TH IEEE International conference on machine
learning and applications (ICMLA) (pp. 635, 642). IEEE}




\bibitem[Yu, W., Kim, I. Y.,  Mechefske, C. (2021)]{Yu:21}
Yu, W., Kim, I. Y., Mechefske, C. ,2021. Analysis of different RNN autoencoder variants for time
series classification and machine prognostics. {Mechanical Systems and Signal Processing, 149, 107322}


\bibitem[Bagnall, A., Lines, J., Bostrom, A., Large, J.,  Keogh, E.(2017)]{Bagnall:17}
Bagnall, A., Lines, J., Bostrom, A., Large, J.,  Keogh, E. ,2017. The great time series classification
{bake off: a review and experimental evaluation of recent algorithmic advances. Data mining and
knowledge discovery, 31(3), 606-660}


\bibitem[Gamboa, J. C. B. (2017)]{Gam:17}
Gamboa, J. C. B. 2017. Deep learning for time series analysis.


\bibitem[A. Q. Khan,S. M. Qureshi (2020)]{A khan:20}
A. Q. Khan,S. M. Qureshi (2020) . Dynamical properties of some rational systems of difference equations.{Mathematical method in the applied science}

\bibitem[José M. Amigó , Michael Small
 (2017)]{Small:17}
A. Q. Khan,S. M. Qureshi (2017) .Mathematical methods in medicine: neuroscience, cardiology and pathology
{National library of medcin}

\bibitem[AL Goldberger,E Goldberger
 (1977)]{Goldberger:77}
AL Goldberger,E Goldberger
 (1977) .Clinical electrocardiac ,Louis MO

\bibitem[Jiapu Pan; Willis J. Tompkins
(1985)]{Jpn:85}
Jiapu Pan; Willis J. Tompkins
(1985).A Real-Time QRS Detection Algorithm


\bibitem[P.E. McSharry , G.D. Clifford, L. Tarassenko, L.A. Smith
(2003)]{P.E:03}
P.E. McSharry , G.D. Clifford, L. Tarassenko, L.A. Smith
(2003).A dynamical model for generating synthetic electrocardiogram signals .{IEEE Transactions on Biomedical Engineering ( Volume: 50, Issue: 3)}


\bibitem[M Malik , AJ camm(1995)]{AJ:95}
M Malik , AJ camm(1995). Heart rate variability ,Armonk,NY :futura

 
 \bibitem[Task force of the european society of cardiology(1996)]{North:96}
Task force of the european society of cardiology(1996). Heart rate variability :standard of measurement physiological interpretation and clinical use {the north americain society of pacing and electrophysiology,sophia antipolis,France}

 
 \bibitem[M H Crawford,S Bernstein and P Deedwania(1999)]{M H Crawford:96}
M H Crawford,S Bernstein and P Deedwania(1999). ACC\ AHA , guidlines for ambulatory ,electrocardiography.{circulation vol 100, pp 886-893}

 
 
 \bibitem[Jianwen LuoKui ,YingLijing BaiLijing Bai(2005)]{Jianwen LuoKui:05}
Jianwen LuoKui ,YingLijing BaiLijing Bai(2005). Savitzky–Golay smoothing and differentiation filter for even number data

 
 
 \bibitem[C. Li, C. Zheng, and C. Tai(1995)]{Tai:95}
 C. Li, C. Zheng, and C. Tai(1995).Detection of ECG characteristic
points using wavelet transforms. {IEEE Transactions on Biomedical Engineering, vol. 42, no. 1, pp. 21–28}
 
 
\bibitem[Z. Tang, G. Zhao, and T. Ouyang(2021)]{Z.Tang} 
Z. Tang, G. Zhao, and T. Ouyang(2021).Two-phase deep learning
model for short-term wind direction forecasting,{Renewable
Energy, vol. 173, pp. 1005–1016}

 
 \bibitem[D. A. Winter, P. M. Rautaharju, and H. K. Wolf(1997)]{Wolf:97}, D. A. Winter, P. M. Rautaharju, and H. K. Wolf(1997).Measurement and characteristics of over-all noise content in exercise
electrocardiograms,{American Heart Journal, vol. 74, no. 3,
pp. 324–331, 1967}.
 
 
 \bibitem[M. Thomas, M. K. Das, and S. Ari(2015)]{Ari:15}
 M. Thomas, M. K. Das, and S. Ari(2015).
 Automatic ECG arrhythmia classification using dual tree complex wavelet based features,{
AEU - International Journal of Electronics and Communications,
vol. 69, no. 4, pp. 715–721, 2015}
 
\bibitem[Yan Fang ,
 Jianshe Shi, Yifeng Huang, Taisheng Zeng,Yuguang Ye ,
Lianta Su, Daxin Zhu,and Jianlong Huang(2022)]{Yan Fang:22}
 Yan Fang ,
 Jianshe Shi, Yifeng Huang, Taisheng Zeng,Yuguang Ye ,
Lianta Su, Daxin Zhu,and Jianlong Huang(2022).Electrocardiogram Signal Classification in the Diagnosis of Heart
Disease Based on RBF Neural Network .{Hindawi
Computational and Mathematical Methods in Medicine
Volume 2022, Article ID 9251225, 9 pages}
 
 \bibitem[Claudia Lainscsek and Terrence J. Sejnowski(2013)]{Claudia Lainscsek:13}
 Claudia Lainscsek and Terrence J. Sejnowski(2013).Electrocardiogram classification using delay differential equations
\doi{doi: 10.1063/1.4811544}
 
 \bibitem[R. P. AGARWAL and Y. M. CHOW(1985)]{R.P.Ag:85}
R. P. AGARWAL and Y. M. CHOW(1985). FINITE-DIFFERENCE METHODS FOR BOUNDARY-VALUE
PROBLEMS OF DIFFERENTIAL EQUATIONS
WITH DEVIATING ARGUMENTS 

 
 
 \end{thebibliography}







\end{document}